\documentclass[12pt]{article}
\usepackage{tikz}
\usepackage{graphicx}
\usepackage{multirow}
\usepackage{amsmath}
\usepackage{indentfirst}
\usepackage{mathrsfs}

\input amssym

\newtheorem{lem}{Lemma}[section]%
\newtheorem{thm}[lem]{Theorem}%
\newtheorem{exam}[lem]{Example}%
\newenvironment{pf}{\medskip\noindent{Proof.\hspace{0.2cm}}}{\hfill \qed \newline \medskip}

\newcommand{\qed}{\hfill \mbox{\raisebox{0.7ex}{\fbox{}}}}

 \def\b{\beta}  \def\d{\delta} 
 \def\s{\sigma}

\def\nd{\mathrel{\bigm|\kern-.7em/}}

\begin{document}

\title	{\Large  Semi-3-abelian 3-group}
	
\author{Xuesong Ma$^{1}$ and Wei Xu$^{2}$\\ School of Mathematical Sciences, Captial Normal University\\ Beijing,  100048,People's Republic of China,	
\thanks{$^1$maxues@cnu.edu.cn\; $^2$ xw\_7314@sina.com}}
	
\maketitle	
	\begin{abstract} In this paper, we give some properties of strong semi-3-abelian 3-groups. And we give  semi-3-abelian 3-group  which any proper subgroups are semi-9-abelian, but the group is not semi-9-abelian. Then there exists one class group which is   semi-3-abelian 3-group but  not strongly semi-3-abelian.
	\end{abstract}
	

{\bf  AMS:20D15}	

	\section{Introduction}
	Let $G$ be a finite p-group. For a positive integer $i\in N$,  $G$ is called a {\sl semi-$p^{i}$-abelian p-group}  provided for all $a,b\in G$ such that $(ab)^{p^{i}}=1$ if and only if $a^{p^{i}}b^{p^{i}}=1$.  And $G$ is called a {\sl strongly semi-p-abelian p- group}  if $G$ is a {\sl semi-$p^{i}$-abelian} for all $i\in N$. $G$ is called a {\sl regular p-group} if $(ab)^{p^{i}}=a^{p^{i}}b^{p^{i}}c_{3}^{p^{i}}\cdots c_{m}^{p^{i}}$ where for all $a,b\in G$ and all $i\in N$,  $c_{j}\in <a,b>', j=3,4,\cdots, m$. Denote $\Omega_{i}(G)=<\{g |g^{p^{i}}=1, g\in G\}>$ and $\mho_{i}(G)=<\{g^{p^{i}}| g\in G\}> $. Define $\pi_{i}: G/\Omega_{i}(G) \rightarrow \mho_{i}(G)$ by $\pi_{i}(\overline{g})=g^{p^{i}}$.  The map $\pi_{i}$ is well defined if and only if $G$ is a semi-$p^{i}$-abelian p-group. In \cite{xu1} MY. Xu  proved  that $G$ was a strongly semi-p-abelian p-group if $\pi_{i}$ was a  bijective map for all $i\in N$  and called $G$ has a {\sl regular power structure} if the $G$ satisfied following three properties.
	
	\begin{enumerate}
		\item[(1)] $\mho_{i}(G)=\{g^{p^{i}}| g\in G\}$;
		\item[(2)] $\Omega_{i}(G)=\{g |g^{p^{i}}=1, g\in G\}$;
		\item [(3)]$\pi_{i}$ is a bijective map.
	\end{enumerate}
	
	In \cite{hall1} P. Hall poved that a regular  p-group had  the  regular power structure and that there existed  metacyclic p-groups with the regular power structure but which are not  regular p-groups. There is a gap between them. It is an interesting problem on determining all irregular p-groups which have the regular power structure. $G$ is called $P_3-group$ if all sections of $G$  satisfy $|G:\Omega_{i}(G)|=|\mho_{i}(G)|$ for all $i\in N$.
	$G$ is called $P_{i}-group$ if all sections of $G$ satisfied property (i) where $i=1,2$. In \cite{Amann1} A. Mann proved that   $P_{3}$-group is  $P_{2}$-group and that $P_2$-group is $P_{1}$-group. In \cite{xu2} MY. Xu proved that a semi-2-abelian group was a strongly semi-2-abelian 2-group.  And  classified  the irregular p-group of  two generators 2-groups with the regular power structure.  In \cite{xu1} MY. Xu gave the problem:  For an odd prime p , are semi-p-abelian p-groups strongly semi-p-abelian?  In this paper, we find the semi-3-abelian 3-groups  which any proper subgroups are semi-9-abelian, but the group is not semi-9-abelian in Exmaple \ref{exa}. Obviously this class groups are not  regular. Thus, the answer of the Xu's problem is negative. We also get some necessary properties of strongly semi-3-abelian 3-groups in Theorem \ref{t2} and \ref{t1},  prove them in section 3,  and give some preliminary lemmas in section 2.

\begin{thm}\label{t2} Let $G$ be a strongly semi-3-abelian 3-group and the exponent  $exp(G)\leq 3^{r}$.  If $G$ is generated by $d$ elements, then $c(G)\leq (r-1)(d+1)+3$
\end{thm}
	
\begin{thm}\label{t1} Let $G$ be a semi-3-abelian but semi-9-abelian 3-group which every proper subgroup is semi-9-abelian, then $G_{7}=1$.
\end{thm}

	Furthermore, we prove the following property for a metabelian semi-p-abelian group.
\begin{thm}\label{t3}
	Meta-abelian semi-p-abelian p-groups are strongly semi-p-abelian.
\end{thm}

	The notation and terminology  are standard in this paper.

\section{Preliminaries}
	
	 A p-group $G$ is called $p-abelian$ if $(ab)^{p}=a^{p}b^{p}$ where for any $a,b\in G$ and $G$ is called inner $semi-p^2-abelian\; group$ if each proper subgroup is a $semi-p^2-abelian\; group$ and $G$ is not $semi-p^2-abelian$.  The following lemmas  are necessary for main results' proof.
	
\begin{lem}\label{ct1}
		Let $G$ be a finite meta-abelian p-group which generated by two elements, then $G$ is p-ableian if and only if $exp(G')\leq p$ and $c(G)<p$.
\end{lem}

There are some basic properties of semi-p-abelian group.

\begin{lem}\label{ct2}
		Let $G$ be a semi-$p^i$-abelian p-group. Then
		\begin{enumerate}
			\item[(1)] each subgroup of $G$ is a semi-$p^i$-abelian p-group;
			\item[(2)]  $[a^{p^i},b]=1$ if and only if  $[a,b]^{p^i}=1$,$\forall a,b\in G$.
		\end{enumerate}
\end{lem}

 Let $G$ be a strongly semi-p-abelian group, obviously  $G$ is semi-$p^2$-abelian. Conversly, we also have following lemma.

	\begin{lem}\label{ct3}
		If any meta-abelian semi-p-abelian group is semi-$p^2$-group,then a meta-abelian semi-p-abelian is strongly semi-p-abelian.
		\end{lem}
	\begin{pf} Let $G$ be a meta-abelian semi-p-abelian group. Suppose any meta-abelian semi-p-abelian groups are semi-$p^2$-abelian. By induction on $|G|$, since $G$ is semi-$p^2$-abelian, then $G/\Omega_{1}(G)$ is a meta-abelian semi-p-abelian group. From assumption, $G/\Omega_{1}(G)$ is a strongly semi-p-abelian.Then $G$ is strongly semi-p-abelian.
	\end{pf}

	Following P.Hall formula is crucial for the research of p-group.

	\begin{lem}\label{ct4}
		Let $G$ be a p-group, then  $\forall a,b\in G$
		$$a^{p}b^{p}=(ab)^pd_{2}^{ {p \choose 2 }}d_{3}^{{p \choose 3 }}\cdots d_{p-1}^{{p \choose p-1 }}d_{p},$$
	\end{lem}where $d_{i}\in H_{i}$ and $H=<a,b>$.

\begin{lem}\label{ct5}
	Let $G$ be a semi-p-abelian and inner semi-$p^2$-abelian group. Then $exp(G')\leq p^2$,$\mho_{2}(G)\leq Z(G)$ and $\mho_{1}(G)'=1$ .
\end{lem}
\begin{pf}  Without loss generality, assume $G=<a,b>$ where  $(ab)^{p^2}=1$ and $a^{p^2}b^{p^2}\neq1$ or $a^{p^2}b^{p^2}=1$ and $(ab)^{p^2}\neq1$.
	
	Suppose $(ab)^{p^2}=1$ and $a^{p^2}b^{p^2}\neq1$,  then $[(ab)^{p^2},b]=1$.
	Notice that $<[ab,b],ab>$ is a proper subgroup of $G$, then $<[ab,b],ab>$ is a semi-$p^2$-abelia p-group. Furthermore $[ab,b]^{p^2}=1$.
	Since $G'$ is a semi-$p^2$-abelian p-group, then $exp(G')\leq p^2$. Thus $\mho_{2}(G)\leq Z(G)$ from Lemma\ref{ct2}. Since $G$ is semi-p-abelian p-group, $\mho_{1}(G)'=1$.
	
	It is simliar for  $a^{p^2}b^{p^2}=1$ and $(ab)^{p^2}\neq1$.

\end{pf}

\section{The proof of main results}

 In order to prove the Theorem \ref{t1} and \ref{t2},  firstly we need following lemmas.

	\begin{lem} \label{tg1}
		Let $G$ be a semi-3-abelian 3-group and $x\in\Omega_{1}(G)$, then  for all $a,b\in G$, $$[x,b,a][x,a,b]=1.$$
	\end{lem}
	
	\begin{pf}
		Let $x\in  \Omega_{1}(G)$, $ y\in G$ and $H=<x,y>$,  then $[x,y,y]=1$ by lemma\ref{ct1}.
		$\forall a,b\in G$, since $[x,a], [x,b]\in  \Omega_{1}(G)$, $[x,a,b,b]=[x,b,a,a]=1$.
		Notice $$[x,ab,ab]=[x,b,b]^{[x,a]^{b}}[x,b,a]^{[x,a]b}[x,a,b]^{[x,a,b]}[x,a,a]^{b[x,a,b]}[x,a,b,a]^{b}.$$
		Since $[x,a,a]=[x,b,b]=1$, then $$[x,ab,ab]=[x,b,a]^{[x,a]b}[x,a,b][x,a,b,a]^{b}.$$
		For $[x,a,b,b]=[x,a,b,[x,a]]=1$, we have
		\begin{equation*}
			[x,b,a][x,a,b][x,a,b,a]^{[a,x]}=1.\tag{3.1.1}\label{s1}
		\end{equation*}
		Similarly,\begin{equation*}
			[x,a,b][x,b,a][x,b,a,b]^{[b,x]}=1.\tag{3.1.2}\label{s2}
		\end{equation*}
		
		From Witt lemma, we have
		\begin{equation*}
			[[x,b,a,b],[b,x]]^{b^{-1}}[b^{-1},[x,b],[x,b,a]]^{[b,x]}[[b,x],[x,b,a]^{-1},b^{-1}]^{[x,b,a]}=1.
		\end{equation*}
		
		By $[x,b,b]=1$ and $[[x,b,a],[x,b]]=1$, then $[b^{-1},[x,b]]=[[b,x],[x,b,a]^{-1}]=1$. So $[[x,b,a,b],[b,x]]=1$.
		
		Similarly, then $[[x,a,b,a],[a,x]]=1$. From \eqref{s1} and \eqref{s2}, then
		\begin{align*}
			&[x,b,a][x,a,b][x,a,b,a]=1,\\
			&[x,a,b][x,b,a][x,b,a,b]=1.
		\end{align*}
		
		Then $[x,a,b,a]=[x,b,a,b]=1$ and $[x,b,a][x,a,b]=1$.

	\end{pf}
	
From lemma\ref{tg1}, $\Omega_{1}(G)\leq Z_{3}(G)$ if $G$ is a semi-3-abelian group which generated by two elements. But there is a general conclusion by lemma\ref{tg1}.

\begin{lem}\label{tg2}
	Let G be a semi-3-abelian group. If $x\in\Omega_{1}(G)$ and $g_{1},g_{2},\cdots,g_{k}\in G$, then $$[x,g_{\s{(1)}},g_{\s{(2)}},\dots,g_{\s{(k-1)}},g_{\s{(k)}}]=[x,g_{1},g_{2},\dots,g_{k-1},g_{k}]^{2^{\d{(\s)}}},$$ where $\s$ is a permutation on t elements and $\d$ is a function
	
	\begin{equation*}
		\d{(\s)}=\left\{\begin{array}{rcl}
			0&&{\s\in Alt_{k}}\\
			1&&{\s\not\in Alt_{k}}
		\end{array}\right.
	\end{equation*}
\end{lem}
\begin{pf}
	 By induction on $k$, we prove \begin{align*}
		[x,g_{1},\cdots,g_{i},g_{i+1},\cdots, g_{k}]&
		=[x,g_{1},g_{2},\cdots,g_{i-1},g_{i+1},g_{i},g_{i+2},\dots,g_{k}]^2.
	\end{align*}
	Notice $[x,b,a][x,a,b]=1$ for $a,b\in G$, then assume  the result is true for  $k-1$.
	
	%%%%%%%%%%%%%%%%%%%%%%%%%%
	
	If $i\leq k-1$, then  \begin{equation*}[x,g_{1},\cdots,g_{i},g_{i+1},\cdots, g_{k-1}]=[x,g_{1},g_{2},\cdots,g_{i-1},g_{i+1},g_{i},g_{i+2},\dots,g_{k-1}]^2.\end{equation*}
	Since $<[x,g_{1},g_{2},\cdots,g_{i-1},g_{i+1},g_{i},g_{i+2},\dots,g_{k-1}],g_{k}>$ is a 3-abelian group generated by two elements, we have \begin{align*}[[x,g_{1},g_{2},\cdots,g_{i-1},g_{i+1},g_{i},g_{i+2},\dots,g_{k-1}]^2,g_{k}]=[x,g_{1},g_{2},\cdots,g_{i-1},g_{i+1},g_{i},g_{i+2},\dots,g_{k-1},g_{k}]^2.\end{align*}
	Then $[x,g_{1},\cdots,g_{i},g_{i+1},\cdots, g_{k}]
	=[x,g_{1},g_{2},\cdots,g_{i-1},g_{i+1},g_{i},g_{i+2},\dots,g_{k}]^2$.
	
	If $i=k-1$,since $[x,g_{1},\cdots,g_{k-2}]\in \Omega_{1}(G)$,then $[x,g_{1},g_{2}\cdots,g_{k-2},g_{k},g_{k-1}]
	=[x,g_{1},g_{2},\cdots,g_{k-2},g_{k-1},g_{k}]^2$.
	
	Now we proved $[x,g_{1},\cdots,g_{i},g_{i+1},\cdots, g_{k}]
	=[x,g_{1},g_{2},\cdots,g_{i-1},g_{i+1},g_{i},g_{i+2},\dots,g_{k}]^2.$ Since $exp(\Omega_{1}(G))\leq 3$, then$$[x,g_{\s{(1)}},g_{\s{(2)}},\dots,g_{\s{(k-1)}},g_{\s{(k)}}]=[x,g_{1},g_{2},\dots,g_{k-1},g_{k}]^{2^{\d{(\s)}}}$$where $\s$ is a permutation of t elements and $\d$ is a function
	\begin{equation*}
		\d{(\s)}=\left\{\begin{array}{rcl}
			0&&{\s\in Alt_{k}}\\
			1&&{\s\not\in Alt_{k}}
		\end{array}\right.
	\end{equation*}
	%%%%%%%%%%%%%%%%%%%
\end{pf}

	\begin{lem}\label{tg3}
		Let $G=<a_{1}, a_{2}, \cdots, a_{n}>$ be a semi-3-abelian group and $|G/\Phi(G)|=3^{n}$, then $\Omega_{1}(G)\leq Z_{n+1}(G)$. Assume  $x=a_{i}^{m}\in \Omega_{1}(G)$ where $m\in N$, then $x\in Z_{n}(G)$.
	\end{lem}
	
	\begin{pf}Let $g_{1},g_{2},\cdots, g_{k}\in G$ and $x\in  \Omega_{1}(G)$ where $k\in N$. From lemma\ref{tg2},we know $[x,g_{1},g_{2},\cdots,g_{i},\cdots,g_{k}]$ is power of $[x,g_{1},g_{2},\cdots,g_{i-1},g_{i+1}\cdots,g_{k},g_{i}]$.Then $[x,g_{1},g_{2},\cdots,g_{k}]$ commutate with $g_{i}$ for $1\leq i\leq k$,  $x\in Z_{n+1}(G)$.

If $x={a_{i}}^{m}$ for $1\leq i\leq n$,then  \begin{equation*}
	[x,a_{\b{(i)}},a_{\b{(i+1)}},a_{\b{(i+2)}},\dots,a_{\b{(n)}},a_{\b{(1)}},\dots,a_{\b{(i-1)}}]=[x,a_{i},a_{i+1},a_{i+2},\dots,a_{n},a_{1},\dots,a_{i-1}]^{2^{\d{(\b)}}}=1
\end{equation*}
where $\b$ is a permutation  of $n$ elements. So $x\in Z_{n}(G)$
	\end{pf}
	
From lemma\ref{tg3},  the nilpotent class of strongly semi-3-abelian group  is  boundded  by the number of generator and its exponent.

 {\noindent\bf Proof of Theorem \ref{t2}}.

\begin{pf}
	Since G is a strongly semi-3-abelian group, then $exp(G/\Omega_{r-1}(G))=3$. Thus $G_{4}\leq\Omega_{r-1}(G)$. By lemma3.1,  $G_{5+d}\leq\Omega_{r-2}(G)$. By induction on $exp(G)$, suppose the result is true for $exp(G)\leq r-1$. Then $G_{(r-2)(d+1)+4}(G)\leq\Omega_{1}(G)$. From lemma3.1, we have $G_{(r-1)(d+1)+4}(G)=1$.
	\end{pf}
	
	From now on, we begin to give some characterizations of semi-3-abelian and inner semi-9-abelian group. Firstly, we give a special form of lemma\ref{ct4}.
	\begin{lem}\label{tg4}
		Let G be a group and $a,b$ elements of G. Then$$a^3b^3=(ab)^3[ab,b]^3[ab,b,ab]^{[ab,b]^2}[ab,b,b].$$
		\end{lem}
	
	\begin{lem}\label{tg5}
		Let  $G$ be a finite semi-3-abelian and inner semi-9-abelian group. Then
		
		(1) $G_{3}\leq \mho_{1}(G)$ and $exp(G_{4})\leq 3$;
		
		(2) $G_{7}=1$.
	\end{lem}
	
\begin{pf}Since the nilpotent class of 3-group with exponent 3 generated by two elements is less than 2.We get $G_{3}\leq \mho_{1}(G)$. Further step,$G_{4}\leq[\mho_{1}(G),G]$.From $exp(G')\leq 9$,we know $exp([\mho_{1}(G),G])\leq3$.Then $exp(G_{4})\leq 3$. By lemma\ref{tg3},  $G_{7}=1$.

\end{pf}
	\begin{thm}\label{tg6}
		Suppose G is a semi-3-abelian and inner semi-9-abelian. Then one of following is true.
		
		(1) If there is $a,b\in G$ satisfy $(ab)^9=1$ and $a^9b^9\neq1$,then $[ab,b,b]^3=[ab,b,b,ab,ab]$ and $[ab,b,b,ab,b]\in Z(G)$.
		
		(2) If there is $a,b\in G$ satisfy $a^9b^9=1$ and $(ab)^9\neq1$,then $$[b,a,b,b,b],[b,a,b,b,a],[a,b,a,a,a],[a,b,a,a,b]\in Z(G)$$ and a equality $$[b,a,b,b,b][b,a,b,b,a][a,b,a,a,a][a,b,a,a,b][a,b,a,b,a,b]=1.$$
	\end{thm}
\begin{pf}
	(1)
	Notice that $$(ab)^3\in Z_{2}(G)$$ and by lemma\ref{tg3} then  $[ab,b,(ab)^3]=1$. Furthermore,  $$[ab,b,ab]^3=1$$ for G is semi-$3$-abelian. Then $$[ab,b,b]^3\neq 1.$$
	
	By Witt lemma,  $$[ab,b,b,ab]^{b^{-1}}[b^{-1},(ab)^{-1},[ab,b]]^{ab}[ab,[ab,b]^{-1},b^{-1}]^{[ab,b]}=1.$$
	Since $[b^{-1},(ab)^{-1}]=[b,ab]^{(bab)^{-1}}$ and$[ab,[ab,b]^{-1}]^{[ab,b]}=[ab,b,ab]$. Then $$[ab,b,b,ab]^{b^{-1}}[[b,ab,(bab)^{-1}],[ab,b]]^{ab}[[ab,b,ab]^{[b,ab]},b^{-1}]^{[ab,b]}=1.$$
	Since $G_{7}=1$ and $[ab,b,ab]\in\Omega_{1}(G)$, then $$[ab,b,b,ab,b][[b,ab,(bab)^{-1}],[ab,b],b][ab,b,ab,b^{-1},b]=1.$$
	So $[ab,b,b,ab,b,ab]=1$. It imply that $[ab,b,b,ab]\in Z_{2}(G)$.
	From lemma\ref{tg4},there is the equality:$$[(ab)^3,b,b]=[ab,b,b]^3[ab,b,ab,[ab,b],b]^2[ab,b,ab,ab,b].$$
	Since $[ab,b,ab],(ab)^3\in\Omega_{1}(G)$,  above equation  is simplified:  $$[ab,b,b]^3[ab,b,ab,ab,b]=1.$$ By Witt lemma, $[ab,b,b,ab]^{b^{-1}}[[b^{-1},(ab)^{-1}],[ab,b]]^{ab}[ab,[b,ab],b^{-1}]^{[ab,b]}=1$. Since $[ab,b,b,ab],[ab,b,ab]\in\Omega_{1}(G)$ and $G_{7}=1$, there is  $$[ab,b,b,ab,ab][[b,ab],[ab,b,b],ab][ab,[b,ab],b^{-1},ab]=1.$$ Then $[h^{-1},g]=[h,g]^{g^{-1}}$ for $g,h\in G$. So $$[ab,[b,ab]]=[ab,b,ab]^{[b,ab]}.$$Furthermore, $$[ab,b,b,ab,ab][[b,ab],[ab,b,b],ab]=[ab,b,ab,b,ab].$$  Then $$[[b,ab],[ab,b,b],ab][[ab,b,b,ab],[b,ab]]=1.$$ By Witt formula, $$[ab,b,b,ab,ab]=[ab,b,ab,b,ab]$$ from $[ab,b,b,ab]\in Z_{2}(G)$. From lemma\ref{tg1},we know:$$[ab,b,ab,b,ab]^2=[ab,b,ab,ab,b].$$
	Then there is a equality $[ab,b,b]^3[ab,b,b,ab,ab]^2=1$.This imply $$([ab,b,b])^3=[ab,b,b,ab,ab].$$

	(2) It's easy to know $a^3b^3\in\Omega_{1}(G)$ from $a^9b^9=1$ and $G$ is semi-$3$-abelian. Since $a^3b^3\in\Omega_{1}(G)$, then$$[b^3,a,b][a^3,b,a]=1,[a^3,b,b]=[b^3,a,a]=1.$$
	
	Then
\begin{align}
&[b^3,a,b]=[b,a,b]^3[b,a,b,[b,a],b]^2[b,a,b,b,b],\tag{3.6.1}\label{h1}
\\&[b^3,a,a]=[b,a,a]^3[b,a,b,[b,a],a]^2[b,a,b,b,a],\tag{3.6.2}\label{h2}
\\&[a^3,b,a]=[a,b,a]^3[a,b,a,[a,b],a]^2[a,b,a,a,a],\tag{3.6.3}\label{h3}
\\&[a^3,b,b]=[a,b,b]^3[a,b,a,[a,b],b]^2[a,b,a,a,b].\tag{3.6.4}\label{h4}
\end{align}
	Since $[b^3,a,b],[b^3,a,a],[a^3,b,b],[a^3,b,a]\in Z(G)$ and $G_{7}=1$,  $$[b,a,b,b,b],[b,a,b,b,a],[a,b,a,a,a],[a,b,a,a,b]\in Z(G).$$
	Using Witt formula, then\begin{align*}
		&[b,a,b,[b,a],b]=[a,b,a,[a,b],a]=1,\\&[b,a,b,[b,a],a]=[b,a,b,a,b,a]^2[b,a,b,b,a,a],\\&[a,b,a,[a,b],b]=[a,b,a,b,a,a]^2[a,b,a,a,b,b].
	\end{align*}
	And since $[a,b,a,a],[b,a,b,b]\in\Omega_{1}(G)$,there are:\begin{align*}
		&[b,a,b,[b,a],a]=[b,a,b,a,b,a]^2,\\&[a,b,a,[a,b],b]=[a,b,a,b,a,b]^2.
	\end{align*}
	
	Through mulitiply \eqref{h1} and \eqref{h2},\eqref{h3},\eqref{h4}, then
	$$[b,a,b,b,b][b,a,b,b,a][a,b,a,a,a][a,b,a,a,b][b,a,b,[b,a],a]^2[a,b,a,[a,b],b]^2=1,$$
	further $$[b,a,b,b,b][b,a,b,b,a][a,b,a,a,a][a,b,a,a,b][a,b,a,b,a,b]=1.$$
\end{pf}

{\noindent\bf Proof of Theorem 1.3}.

\begin{pf}Suppose $G$ is a meta-abelian semi-p-abelian  and inner semi-$p^2$-abelian group. Then  $\forall a,b\in G$ such that $(ab)^{p^2}=1$ but $a^{p^2}b^{p^2}\neq1$ or $a^{p^2}b^{p^2}=1$ but $(ab)^{p^2}\neq1$. $G=<a, b>$. Since $G$ is a meta-abelian semi-p-abelian group generated by two elements,then $G_{p}\leq\mho_{1}(G)$. Furthermore, since $G_{p+1}\leq[\mho_{1}(G),G]$, then $exp(G_{p+1})\leq p$. Then $G_{2p}=1$.
	
	 Case(1). Since $(ab)^p\leq\Omega_{1}(G)$, then we get some equalities
	\begin{equation*}
		[(ab)^{p},b,\underbrace{(ab)^{s}b,\cdots,(ab)^{s}b}_{p-2}]=1
	\end{equation*}where $1\leq s\leq p-1$.
	Denote  $x\equiv y$ if $x^{-1}y\in G_{p+1}^{p}G_{2p}$ for $x,y \in G$.
	
	Since $[x^{p},x_{1},x_{2},\cdots,x_{p}]\equiv[x,x_{1},x_{2},\cdots,x_{p}]^{p}$, then we get $$[(ab)^{p},\underbrace{b,\cdots,b}_{i_{1}},\underbrace{ab,\cdots,ab}_{i_{2}}]=1$$
	by Vandermonde Matrix on field of $p$ elements, where $i_{1}+i_{2}=p-1$.
	Then \begin{equation*}
		[ab,\underbrace{b,\cdots,b}_{i_{1}},\underbrace{ab,\cdots,ab}_{i_{2}}]^p=1\end{equation*}
	where $i_{2}\geq 1$. From lemma\ref{ct4},we have $$[ab,\underbrace{b,\cdots,b}_{p-1}]\neq1.$$
	We can notice $[(ab)^{p},\underbrace{b,\cdots,b}_{p-1}]=1$. Since $$[(ab)^{p},\underbrace{b,\cdots,b}_{p-1}]\equiv [(ab),\underbrace{b,\cdots,b}_{p-1}]^{p}[ab,\underbrace{b,\cdots,b}_{p-1},\underbrace{ab,\cdots,ab}_{p-1}]$$then $$[(ab),\underbrace{b,\cdots,b}_{p-1}]^{p}[ab,\underbrace{b,\cdots,b}_{p-1},\underbrace{ab,\cdots,ab}_{p-1}]=1$$This imply $([ab,\underbrace{b,\cdots,b}_{p-1}]ab)^{p}=(ab)^{p}$.Then $[ab,\underbrace{b,\cdots,b}_{p-1}]\in\Omega_{1}(G)$ from $G$ is semi-p-abelian.So $G_{p}^{p}=1$. It's contradiction.
	
	Case(2). Since $G$ is a semi-p-abelian group, then $a^{p}b^{p}\in\Omega_{1}(G)$. From lemma\ref{ct1}, we have $$[a^{p}b^{p},\underbrace{a^{s}b,\cdots,a^{s}b}_{p-1}]=1.$$ By using Vandermonde matrix on field of $p$ elements,we have $$[b^p,a,\underbrace{b,\cdots,b}_{i},\underbrace{a,\cdots,a}_{j}]^{{p-2 \choose j}}[a^p,b,\underbrace{b,\cdots,b}_{i-1},\underbrace{a,\cdots,a}_{j+1}]^{p-2 \choose j+1}=1$$where $i+j=p-2$ and $i\geq1$. We can notice \begin{align*}
		&[a^{p},b,\underbrace{b,\cdots,b}_{i},\underbrace{a,\cdots,a}_{j}]\equiv[a,b,\underbrace{b,\cdots,b}_{i},\underbrace{a,\cdots,a}_{j}]^{p}[a,b,\underbrace{b,\cdots,b}_{i},\underbrace{a,\cdots,a}_{j+p-1}]\\&[b^{p},a,\underbrace{a,\cdots,a}_{j},\underbrace{b,\cdots,b}_{i}]\equiv[b,a,\underbrace{a,\cdots,a}_{j},\underbrace{b,\cdots,b}_{i}]^{p}[b,a,\underbrace{a,\cdots,a}_{j},\underbrace{b,\cdots,b}_{i+p-1}].
	\end{align*}
	Since $[a^{-1},b]=[b,a]^{a^{-1}}$, then \begin{align*}
		&[b,a,\underbrace{a,\cdots,a}_{j},\underbrace{b,\cdots,b}_{i}]^{p}=[a,b,\underbrace{a,\cdots,a}_{j},\underbrace{b,\cdots,b}_{i}]^{-p}\\&[b,a,\underbrace{a,\cdots,a}_{j},\underbrace{b,\cdots,b}_{i+p-1}]=[a,b,\underbrace{a,\cdots,a}_{j},\underbrace{b,\cdots,b}_{i+p-1}]^{-1}.
	\end{align*}
	
	We get $[a^{p},b,\underbrace{b,\cdots,b}_{i},\underbrace{a,\cdots,a}_{j}]=[b^{p},a,\underbrace{b,\cdots,b}_{i},\underbrace{a,\cdots,a}_{j}]^{-1}$. Then $$[a^p,b,\underbrace{b,\cdots,b}_{i},\underbrace{a,\cdots,a}_{j}]^{p-2\choose j}=[a^p,b,\underbrace{b,\cdots,b}_{i-1},\underbrace{a,\cdots,a}_{j+1}]^{p-2\choose j+1}.$$ Since $[a^p,b,\underbrace{a,\cdots,a}_{p-2}]=[b^p,a,\underbrace{a,\cdots,a}_{p-2}]=1$, then $[a^p,b,\underbrace{b,\cdots,b}_{i},\underbrace{a,\cdots,a}_{j}]=1$. This imply $a^{p},b^{p}\in Z_{p-1}(G)$.Then $G_{p}^{p}=1$.
	It's contradiction.
	
	From proof above, any meta-abelian semi-p-abelian groups are semi-$p^2$-abelian. By lemma\ref{ct3}, any meta-abelian semi-p-abelian groups are strongly semi-p-abelian.
	
\end{pf}

 Let $G$ be a semi-3-abelian and inner semi-9-abelian group, by lemma\ref{ct2} (2), then $exp(G')\leq 9$. Since lemma\ref{tg2} and \ref{tg3},$exp(G_3)=9$. Take $ G=<a_1,a_2>$ and $[a_1,a_2]^9=1$. By theroem\ref{tg6} (1) , the order of $a_1$, $o(a_1)=9$ and $o(a_2)=3^n$ where $n\geq 3$. Then $o([a_1,a_2,a_1])=3$ and  $o([a_1,a_2,a_2])=9$ . So $[a_1,a_2,a_1,a_i,a_i]=1$ for $i=1,2$ by lemma\ref{ct1}. By theroem\ref{tg6} (1), $[a_1,a_2,a_2]^3=[a_1,a_2,a_2,a_1,a_1]$ . Since theroem\ref{tg5} and \ref{t3}, $4<c(G)\leq 6$. When $c(G)=5$, $[a_1,a_2]^3\notin G_4$ and $G_4$ is a  elementary 3-group with rank 9, we can get Example\ref{exa}. 

	\begin{exam}\label{exa}
			
		 Let $G=<a_{1},a_{2}|{a_{1}}^{3^2}={a_{2}}^{3^n}=b^{3^2}={c_{1}}^{3}={c_{2}}^{3^2}={d_{r}}^3={e_{v}}^3=1,b=[a_{1},a_{2}],c_{i}=[b,a_{i}],d_{2i-2+j}=[c_{i},a_{j}],e_{2r-2+j}=[d_{r},a_{j}],[e_{v},a_{j}]=1,{c_{2}}^3={e_{5}},e_2e_3=1,e_{1}=e_{4}=1>$, then $G$ is a semi-3-abelian and inner semi-9-abelian 3-group.
		\end{exam}
	
Remark: 	Let $G$ be a semi-3-abelian and inner semi-9-abelian 3-group, then $c(G)=5$ or $6$. Classifying  the semi-3-abelian and inner semi-9-abelian 3-group is next step work.

\end{document}